\renewcommand{\to}{\rightarrow}

\documentclass[11pt]{article}

\usepackage{amsfonts,amssymb}
\usepackage[all, knot]{xy}

\newcommand{\lt}{\triangleleft}

\newcommand{\Z}{\mathbb{Z}}
\newcommand{\R}{\mathbb{R}}

\newcommand{\F}{\mathbb{F}}

\newtheorem{theorem}{Theorem}

\newtheorem{lemma}{Lemma}

\input{epsf.sty}
\usepackage{graphicx}

\begin{document}

\title{A knotted $2$-dimensional foam with non-trivial cocycle invariant}

\author{
J. Scott Carter\footnote{Supported by Brain Pool trust}
\\ University of South Alabama
\and Atsushi Ishii \footnote{Supported in part by JSPS KAKENHI Grant Number 24740037}\\ University of Tsukuba}

\maketitle

\begin{abstract} 
By $2$-twist-spinning the knotted graph that represents the knotted handlebody $5_2$, we obtain a knotted foam in $4$-dimensional space with a non-trivial quandle cocycle invariant. 

\end{abstract} 

\section{Introduction}
Knotted foams are to knotted spheres as knotted trivalent graphs are to classical knots. Consider the spine of the tetrahedron that is obtained by embedding four copies of the topological space that is homeomorphic to the alpha-numeric character ${\sc Y}$ in each of the triangular faces and coning the result to the barycenter of the simplex. This two dimensional space (illustrated below the current paragraph), $Y^2$, has a single vertex, four edges, and six $2$-dimensional faces. Three faces are incident to each edge, and a neighborhood of a point in an open edge is homeomorphic to the foam $Y^1 \times [-1,1].$ A {\it $2$-dimensional foam ($2$-foam)} is a compact topological space, $F$, such that any point has a neighborhood that is homeomorphic to a neighborhood of a point in $Y^2$. Thus a foam is stratified into isolated singular points, $1$-dimensional edges at which three sheets meet, and $2$-dimensional faces. The boundary of a foam is a trivalent graph. 
A {\it closed foam} has empty boundary.
Analogous concepts exist in all dimensions. Just as a trivalent graph can be embedded and knotted in $3$-space, a $2$-foam can be embedded and knotted in $4$-dimensional space. 

\begin{center}
\includegraphics[width=2in]{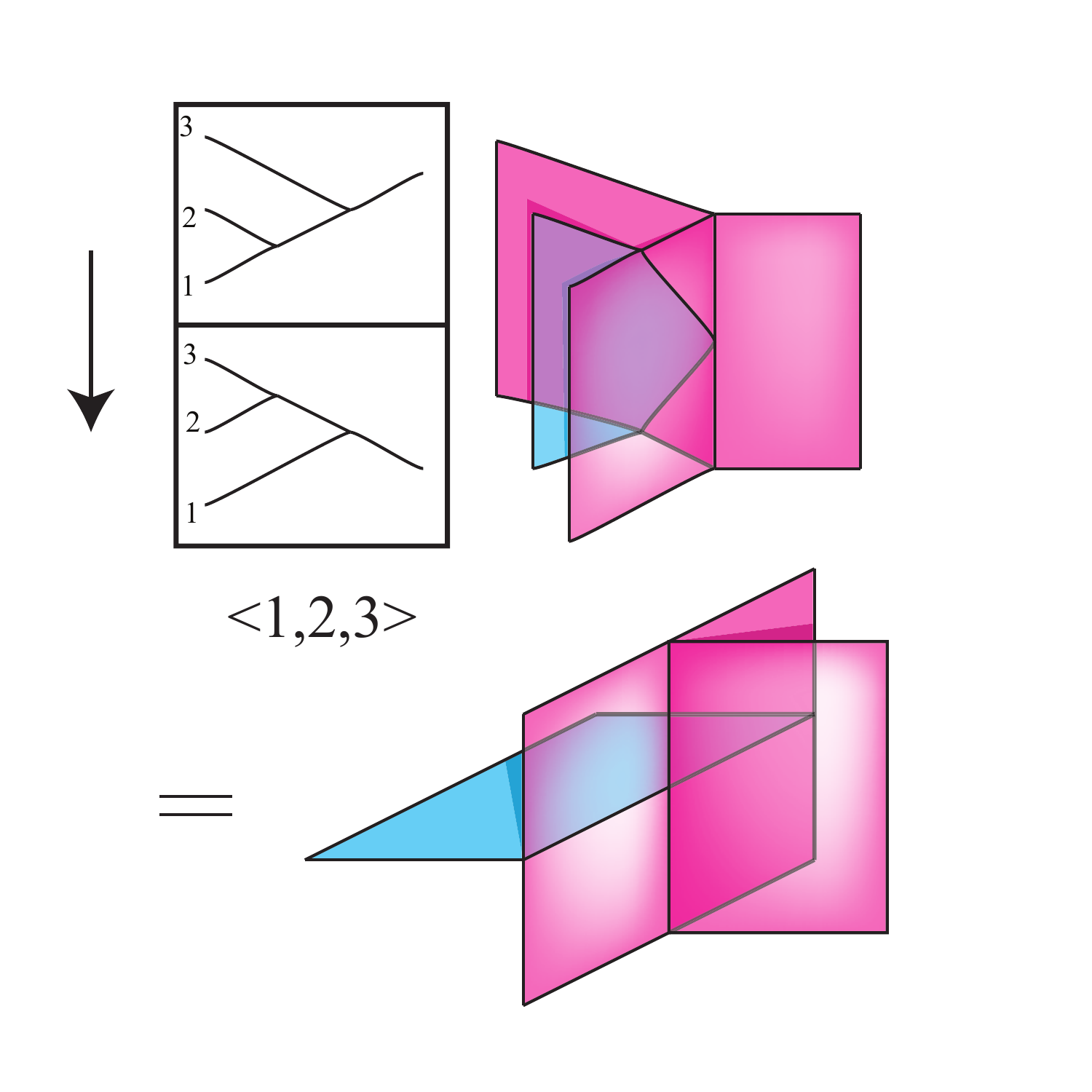}
\end{center}

The space $Y^2$  can be interpreted as a movie of the associativity rule when this is expressed in terms of binary trees. The arrow in the movie presentation indicates a direction determined by the movie that will coincide with sign conventions for the boundary.

An obvious method of constructing examples of knotted $2$-foams is by the method of twist spinning. This operation is achieved by the process that follows that given in \cite{SS} and that is illustrated schematically as follows: 
\begin{center}
\includegraphics[width=3.5in]{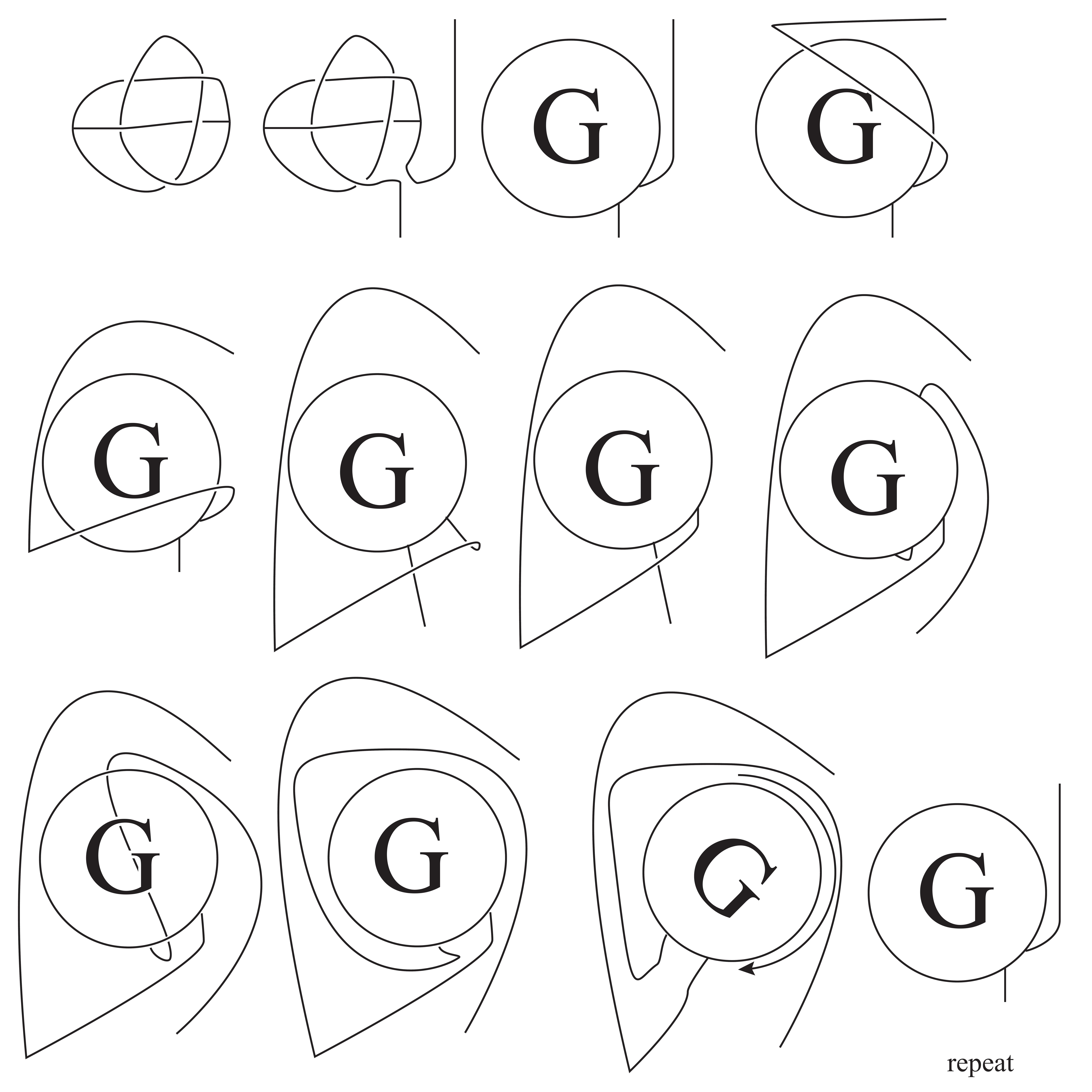}
\end{center} 

The top and bottom edges on the right of the illustration can be capped-off by disks and in this way {\it twist-spinning induces an embedding of a closed foam in $\R^4$.}

Quandle cocycle invariants can be defined for knotted $2$-foams in analogy to the quandle cocycle invariants for knotted trivalent graphs. Here we outline the process in the case that the quandle is an associated quandle to a $G$-family of quandles. 

\subsection*{Acknowledgements} Much of the work for this paper was done in consultation with Masahico Saito who, for reasons of time constraints, was not able to help with the preparation of the manuscript. We anticipate a joint manuscript with Saito-san that more fully develops many of the ideas herein. 
We are also grateful for conversations with Yongju Bae, Seiichi Kamada, Kanako Oshiro, and Shin Satoh as well as the students at the TAPU workshops. This pa?per was studied with the support of the Ministry of Education Science and Technology (MEST) and the Korean Federation of Science and Technology Societies (KOFST).

\section{Group families of quandles}

For the idea of a $G$-family of quandles, we follow the presentation in \cite{IIJO}. Let $G$ denote a group, and let $X$ denote a set 
upon which there is a family of binary operations $\lt_g: X\times X \rightarrow X$ --- one for each element $g \in G$ such that the following properties hold:

\begin{itemize}\item for each $a \in X$ and for each $g \in G$, we have $a \lt_g a =a$;
\item for each $a,b \in X$, and for every $g,h \in G$, we have 
$(a\lt_g b)\lt_h b = a \lt _{gh} b$;
\item the identity element $1 \in G$ induces the trivial operation: $a\lt_1 a =a$;
\item for any $a,b,c \in X$ and for any $g,h \in G$, 
we have $(a\lt_g b) \lt_h c = (a \lt_h c) \lt_{h^{-1}gh} (b \lt_h c).$
\end{itemize}
We read the expression $a \lt_g b$ as, ``$a$ is acted upon by b via the element $g$." The second and third axioms imply that each $\lt_g$ has a left inverse. That is given $g \in G$ and $a, b \in X$, there is a unique $c \in X$ such that $c \lt_g b =a$. To see this let $c = a \lt_{g^{-1}} b$. For fixed $g\in G$, the set $X$ with binary operation $\lt_g$ is a {\it quandle:} every element $a\in A$ is idempotent, the operation is left-invertible, and self distributive. See \cite{JSCQ} for more about quandles.

Given a $G$-family of quandles $\{(X, \lt_g): g\in G \}$, we can define a quandle structure on $X \times G$ via the operation $(a,g)\lt(b,h)=(a \lt_h b, h^{-1}gh)$ where $a,b \in X$ and $g,h\in G$.  This is called {\it the associated quandle} of the $G$-family.

Let $V$ denote a vector space, and let $G$ denote a subgroup of ${\mbox{\rm GL}}(V)$. Then $\{ (V, \lt_M): M\in {\mbox{\rm GL}}(V) \}$ is a $G$-family of quandles under the operations 
$\vec{a}\lt_M \vec{b}= \vec{a}M + \vec{b} - \vec{b}M$, for $\vec{a}, \vec{b} \in V$ and $M \in {\mbox{\rm GL}}(V)$.  (Here we are thinking of elements of $V$ as row vectors.) Note that this idea formalizes the idea of different specializations of the variable $t$ in the definition of the Alexander quandle $a\lt_t b= ta +(1-t)b$. The case that we consider here is $V=\F_3$ with ${\mbox{\rm GL}}(\F_3)=\{\pm 1\}$. It is more convenient to indicate the multiplicative group as $\Z_2 =\{0,1\}$. The quandle operations are $a \lt_0 b=a$ and $a\lt_1 b= 2b-a$. We will denote the associated quandle  $\tilde{R}$

\section{Coloring embedded foams by $\tilde{R}$}

Let $F$ denote a closed embedded foam in $\R^4$.  The elements of $\tilde{R}$ will be indicated as $(a,0)$, $(a,1)$, $(b,0)$, and so forth.  At an edge of $F$ three sheets are coincident. A neighborhood of the edge of $F$ is homeomorphic to $Y \times (-1,1)$. Thus we will refer to the {\it three branches of the foam at an edge}. 
We define a coloring of $F$ by $\tilde{R}$ to be  a function from the set of $2$-dimensional regions of $F$ into the underlying set of $\tilde{R}$ such that: (1) each region is transversely oriented; (2) when three branches at an edge are coincident, then (a) the first components of $\tilde{R}$ are the same (say the color is $a$) and (b) either all the second components are $0$ or exactly one is $0$ and the normal orientations on the sheets colored by $(a,1)$ are consistent. The conditions are indicated below. The double arrows indicate that the normal orientations on the sheets labeled $(a,0)$ can be chosen at will. Also depicted are the possible coincidences of colors at a vertex. In this case the normal directions are also chosen to be consistent along sheets that are colored by $1$. 

\rule{2cm}{0cm}\includegraphics[width=3in]{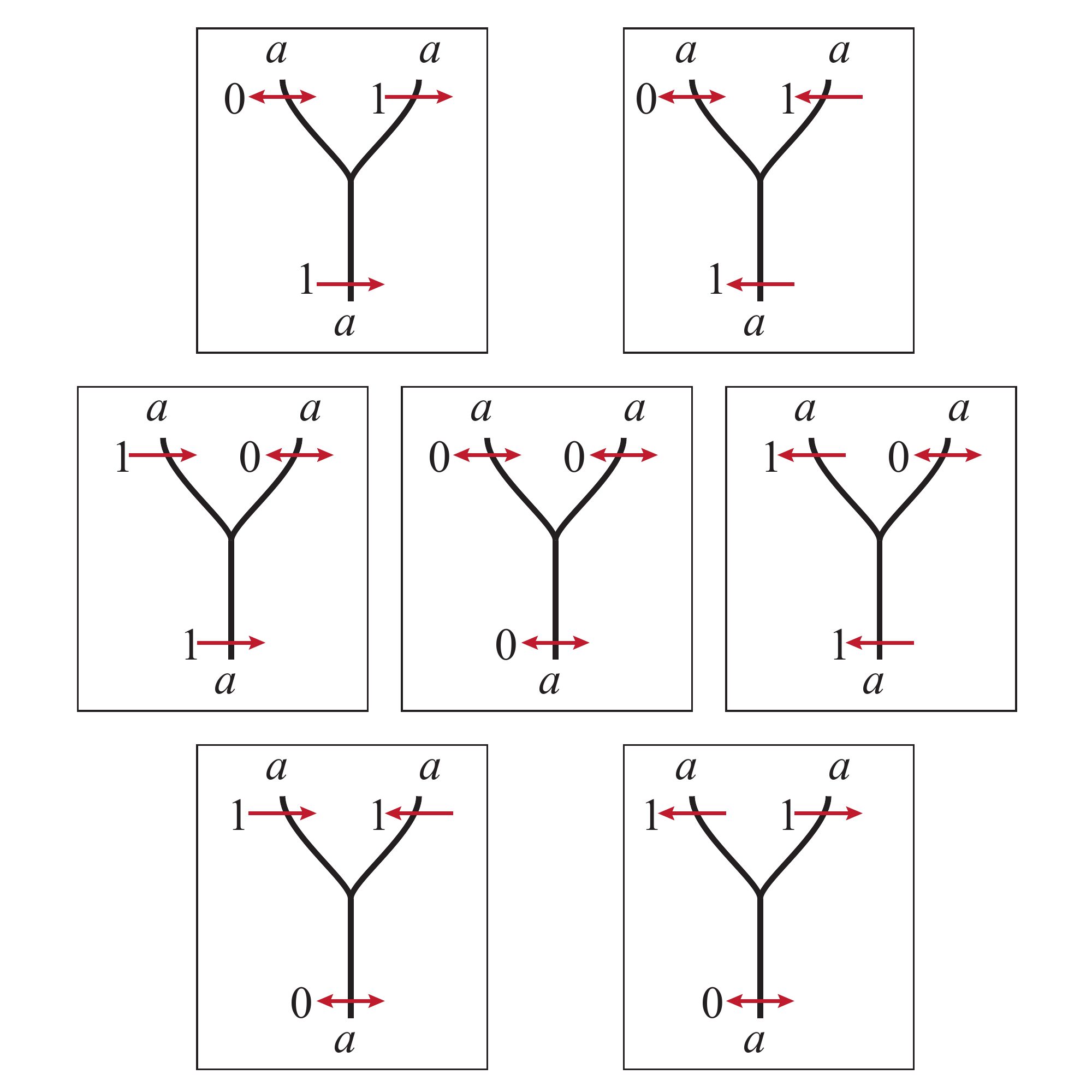}
\includegraphics[width=2.25in]{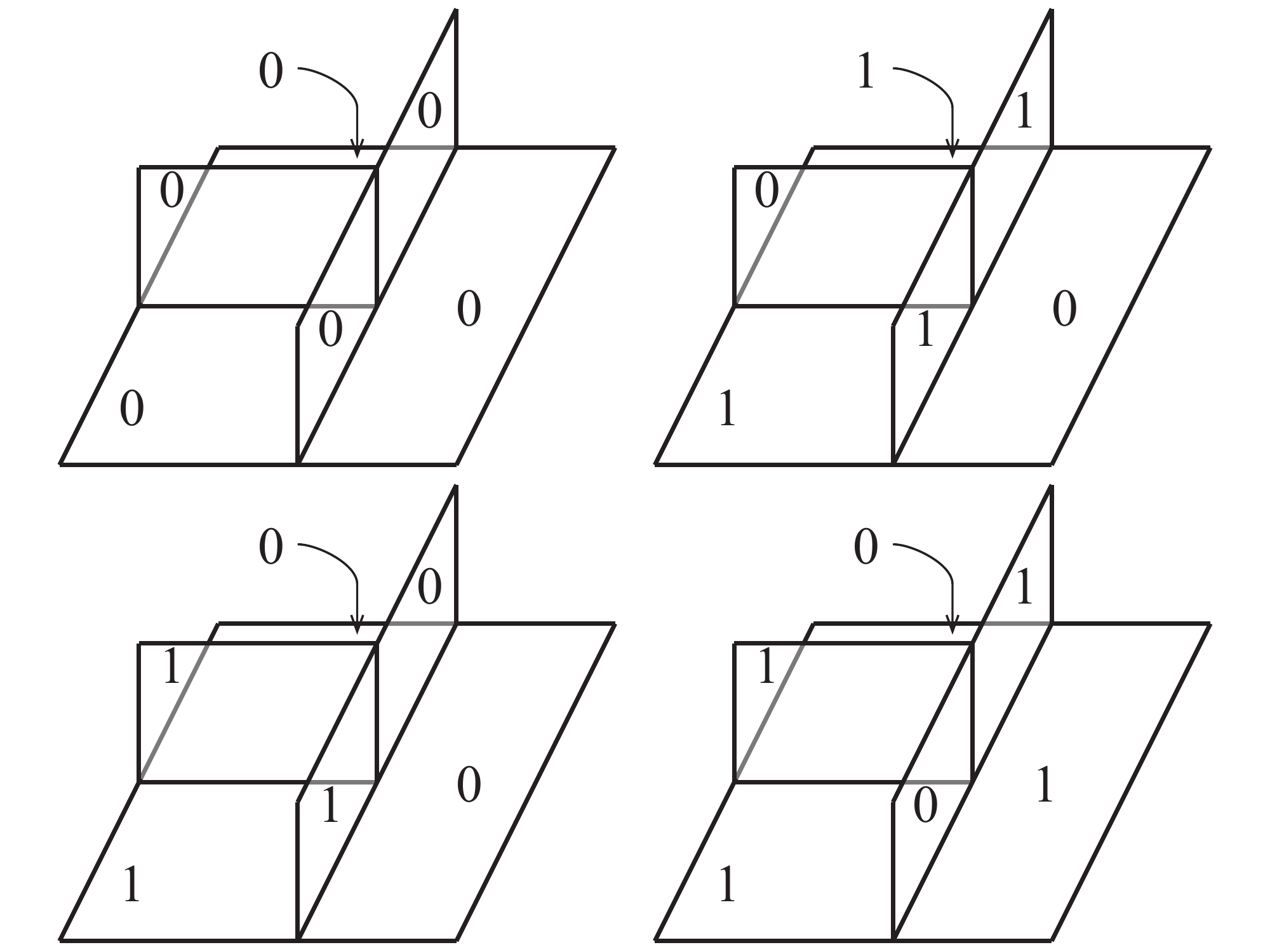}

\section{Homology of $G$-families of Quandles}
When $X$ is a $G$-family of quandles, we define, for each $a \in X$ chain groups, $C_k(a) \{j\}$ --- {\it the set of $k$-chains at $a$ that are off-set by $j$} --- to be the free abelian group generated by $k$-tuples of the form 
$((a, g_{j+1}),(a,g_{j+2}), \ldots, (a, g_{j+k}))$. The element $a$ will be understood in context and to simplify notation, such a chain will be written as $\langle j+1, j+2, \ldots, j+k \rangle$. The associated quandle acts upon chains by 
$$\langle 1, \ldots, k_1\rangle\langle k_1+1, \ldots, k_1+k_2 \rangle \cdots \langle \sum_{i=1}^{\ell-1} k_{i}+1 , \ldots, \sum_{i=1}^{\ell } k_{i} \rangle \lt(j+1)$$
$$ = \langle 1\lt(j+1), \ldots, k_1\lt(j+1)\rangle
\cdots \langle \left(\sum_{i=1}^{\ell-1} k_{i}+1\right) \lt(j+1) , \ldots,\left(\sum_{i=1}^{\ell } k_{i} \right)\lt(j+1) \rangle$$
where $\left(\sum_{i=1}^{\ell } k_{i}\right)=j,$
this indicates the subscripted quantity
$\left( a_\ell, g_{\sum_{i=1}^{\ell } k_{i}}\right)$, and the action is determined by
$(a,g)\lt  (b, h)=(a\lt_h b, h^{-1}g h )$.
The quandle action, then, extends over juxtaposition. The boundary of a chain $\langle j+1, j+2, \ldots, j+k \rangle \in C_k(a) \{j\}$ is computed as follows:
\begin{eqnarray*} \partial \langle j+1, j+2, \ldots, j+k \rangle   &=& \lt(j+1)\langle  j+2, \ldots, j+k \rangle \\ & & + \sum_{\ell=1}^{k-1}  (-1)^\ell\langle j+1,\ldots,  (j+\ell)\cdot(j+\ell+1), \ldots , j+k \rangle \\ & & + (-1)^k \langle j+1, \ldots, j+k-1 \rangle. \end{eqnarray*}
The notation $(j+\ell)\cdot(j+\ell+1)$ indicates the fibre-wise product $(a,g_{j+\ell})\cdot (a, g_{j+\ell+1})= 
(a,g_{j+\ell}\cdot g_{j+\ell+1})$ that is induced by the group structure in $G$. We compute the boundaries under juxtaposition by 
$$\partial(PQ) = (\partial P)Q +(-1)^{{\mbox{\rm dim}} P} P (\partial Q).$$ In general, an $n$-chain is an element of 
$$
C_n= {\displaystyle
{\bigoplus_{(a_1,\ldots, a_\ell)\in X^{\ell}\setminus D}}}
C_{k_1}(a_1)\{0\} \oplus C_{k_2}(a_2)\{k_1\}  \oplus C_{k_\ell}(a_\ell)\left\{\sum_{i=1}^{\ell-1} k_{i}\right\}$$ where the subset $D$ consists of the $\ell$-tuples for which $a_i=a_{i+1}$ for some $i=1,\ldots, \ell-1.$ Here $\sum_{i=1}^{\ell} k_{i}=n$

As usual, a chain, $c$,  is {\it a cycle} if $\partial(c)=0$, and {\it a boundary} if $c=\partial(c')$ for some $c' \in C_{n+1}$. That this defines an homology theory is straight-forward to check and depends upon the associativity in $G$ and upon the self-distributivity of the quandle $X\times G$. 

We will be interested in functions $\alpha$, $\gamma_1$, $\gamma_2$, and $\theta$ that vanish upon the boundaries of certain $4$-cycles. First, we compute the boundaries of generating $3$- and $4$-chains.

 For the generating $3$-chains, we have the following:
 \begin{eqnarray*}
 \partial (\langle 1,2,3 \rangle ) &=& \langle 2,3 \rangle  - \langle 1 \cdot 2, 3 \rangle + \langle 1, 2 \cdot 3 \rangle - \langle 1,2\rangle; \\
  \partial (\langle 1,2\rangle \langle3 \rangle ) &=& 
    \langle 2 \rangle \langle 3 \rangle  - \langle 1 \cdot 2 \rangle \langle 3 \rangle + \langle  1 \rangle \langle 3 \rangle; \\
  \partial (\langle 1\rangle \langle 2,3 \rangle  ) &=& 
     \langle 2,3 \rangle - \langle 2,3 \rangle -\langle 1\lt 2 \rangle \langle 3 \rangle 
    + \langle 1\rangle \langle 2\cdot 3 \rangle -
    \langle 1 \rangle \langle 2 \rangle ; \\
      \partial (\langle 1\rangle \langle 2 \rangle \langle 3 \rangle ) &=& 
     \langle 2\rangle \langle 3 \rangle - \langle 2\rangle \langle3 \rangle -\langle 1\lt 2 \rangle \langle 3 \rangle 
    + \langle 1\rangle \langle  3 \rangle +
    \langle 1\lt 3 \rangle \langle 2 \lt 3 \rangle -  \langle 1\rangle \langle 2 \rangle . 
    \end{eqnarray*}

The three chains listed correspond to the movies to graphs that are illustrated below.
\begin{center}
\includegraphics[width=5in]{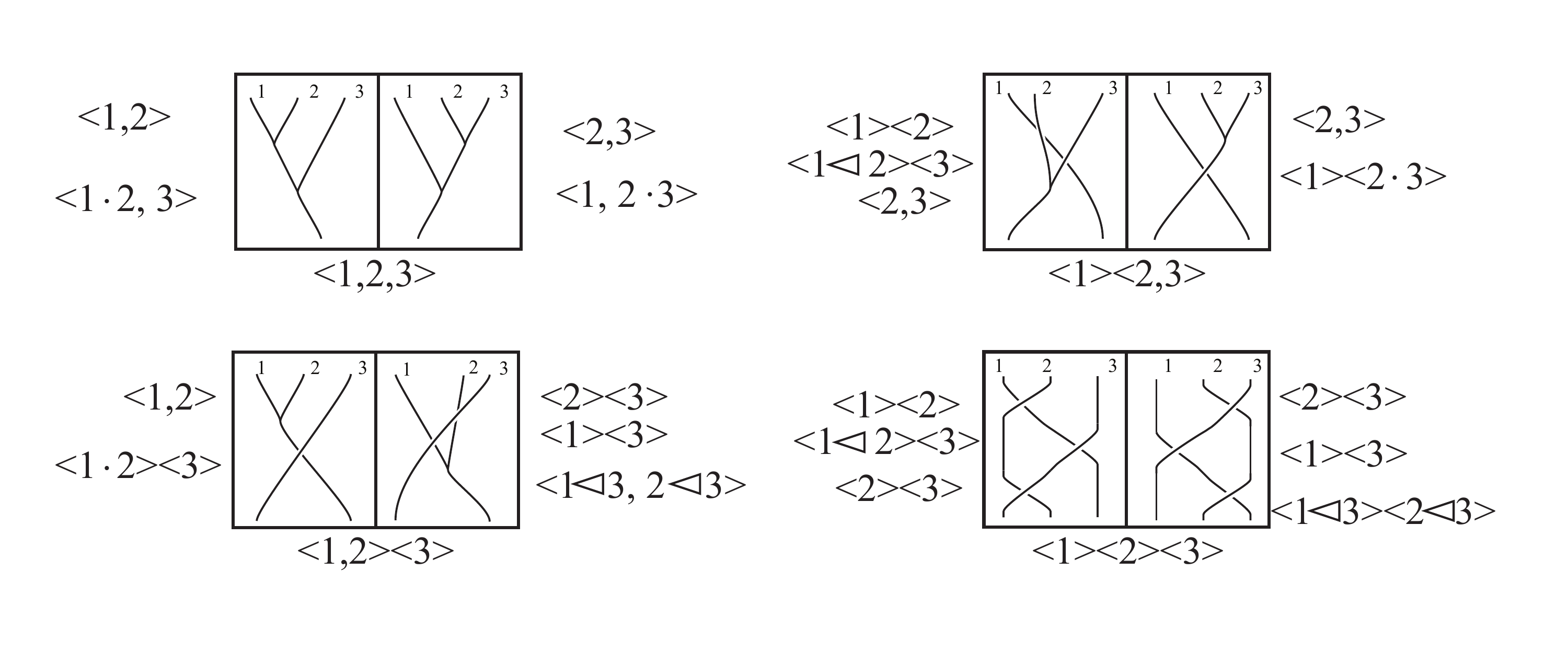}
\end{center}
These are illustrated in broken surface diagram form as follows:

\includegraphics[width=1.75in]{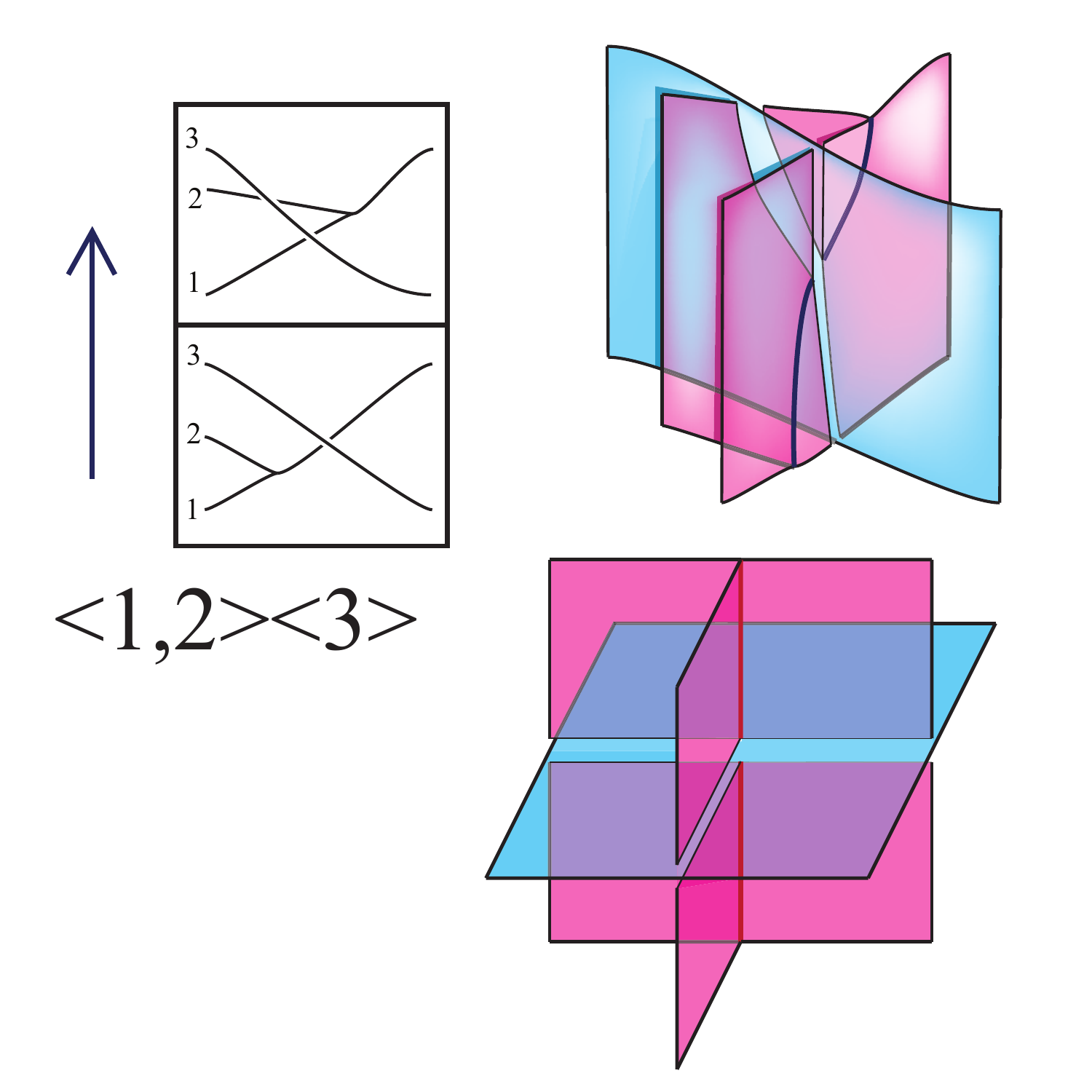}
\includegraphics[width=1.75in]{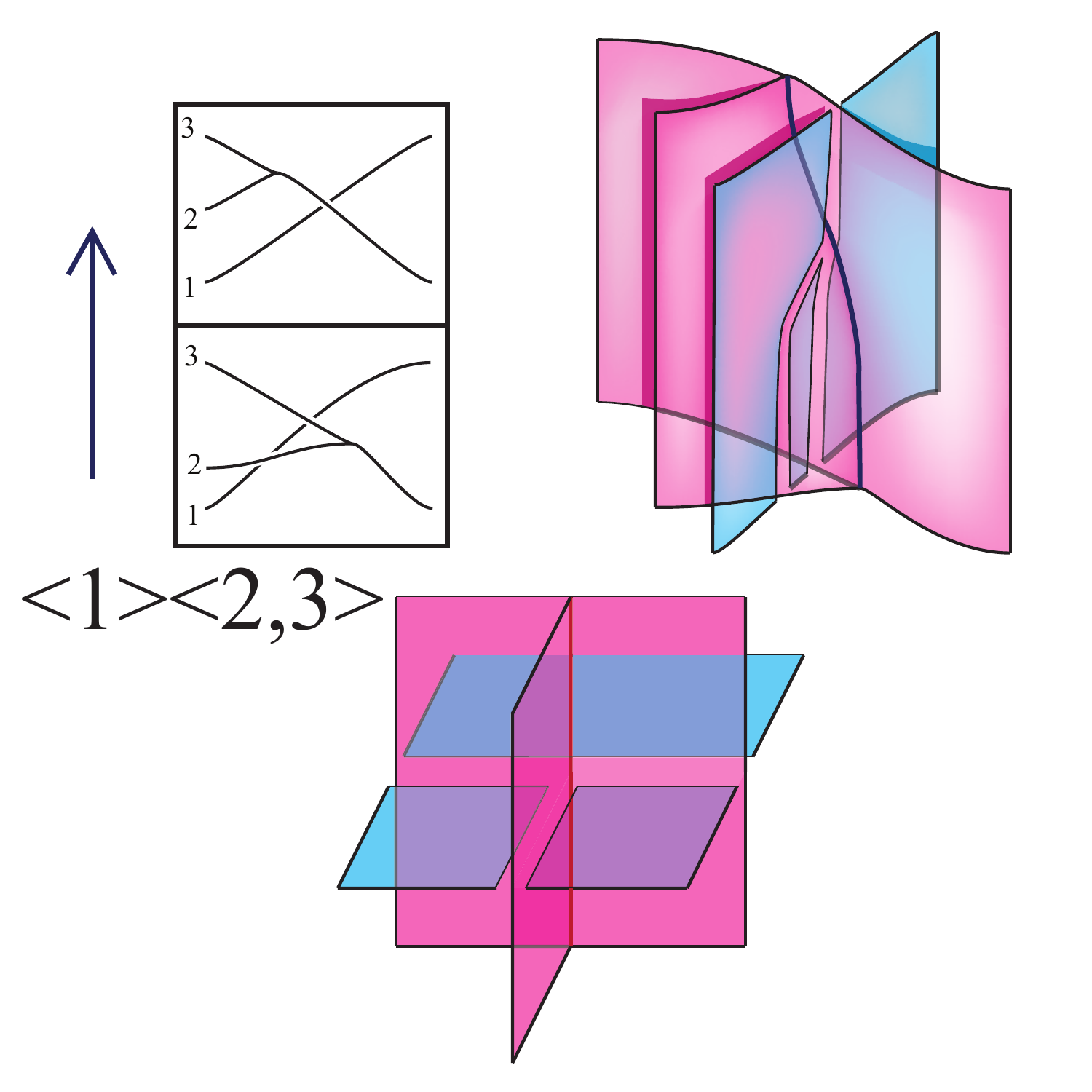}
\includegraphics[width=1.75in]{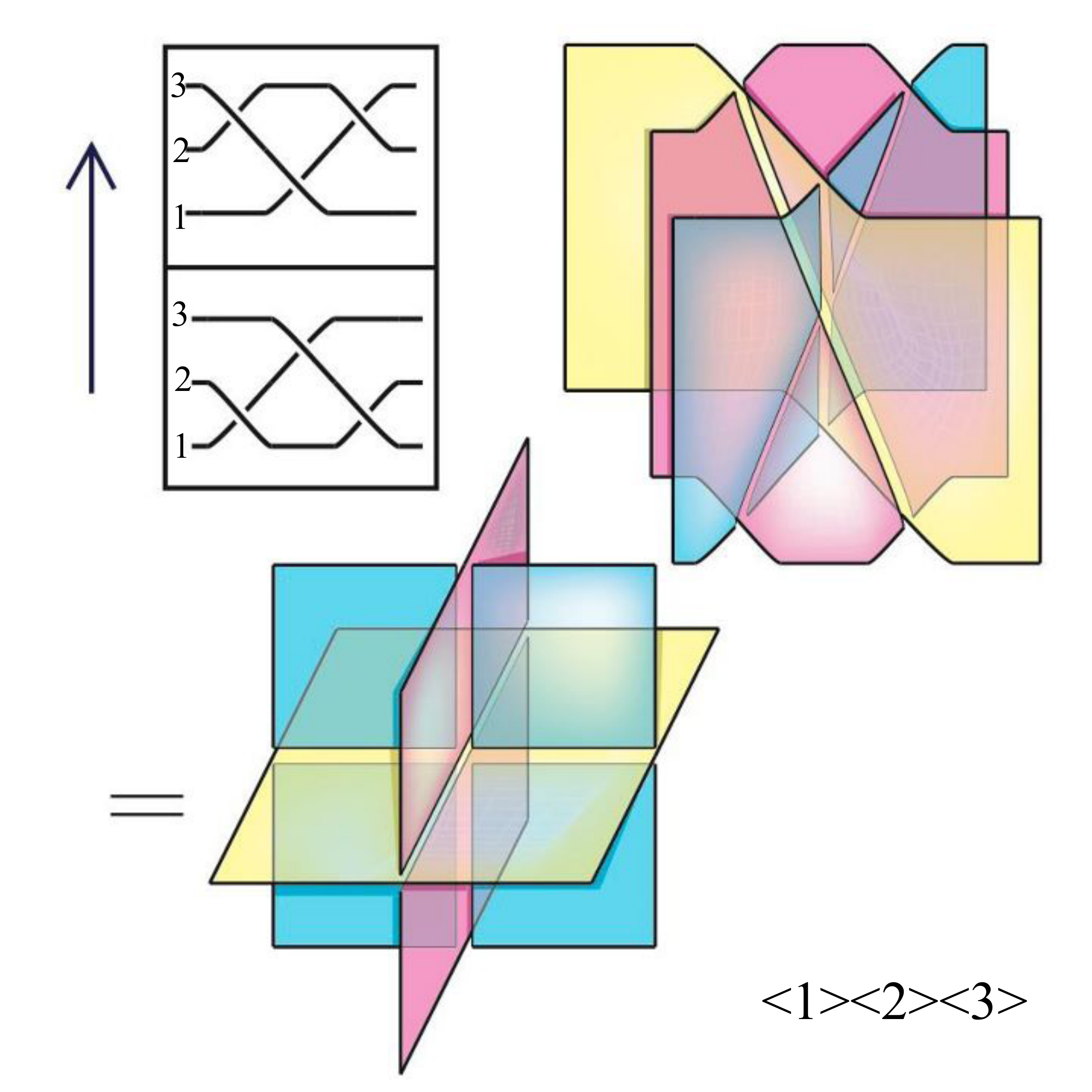}

    \newpage 
    
Meanwhile, for the generating $4$-chains, we have   the following:
   \begin{eqnarray*}
 \partial (\langle 1,2,3,4 \rangle ) &=& \langle 2,3,4 \rangle  - \langle 1 \cdot 2, 3,4 \rangle + \langle 1, 2 \cdot 3,4 \rangle - \langle 1,2, 3 \cdot 4 \rangle + \langle 1,2,3 \rangle ; \\
  \partial (\langle 1,2,3\rangle \langle 4 \rangle ) &=& \langle 2,3 \rangle \langle 4 \rangle - \langle 1 \cdot 2, 3 \rangle \langle 4 \rangle + \langle 1, 2 \cdot 3 \rangle \langle 4 \rangle - \langle 1,2\rangle \langle 4 \rangle 
 - \langle 1\lt 4 , 2 \lt 4 , 3 \lt 4 \rangle + \langle 1 , 2  , 3 \rangle ; \\
   \partial (\langle 1,2\rangle \langle3, 4 \rangle ) &=&
   \langle 2 \rangle \langle 3,4 \rangle  - \langle 1 \cdot 2 \rangle \langle 3,4 \rangle + \langle  1 \rangle \langle 3,4 \rangle 
  + \langle 1\lt 3,2\lt3 \rangle \langle 4 \rangle -
  \langle 1,2 \rangle \langle 3 \cdot 4 \rangle + \langle 1,2 \rangle \langle 3  \rangle ; \\
    \partial (\langle 1,2\rangle \langle3 \rangle \langle 4 \rangle ) &=&
   \langle 2 \rangle \langle 3 \rangle \langle 4 \rangle  - \langle 1 \cdot 2 \rangle \langle 3 \rangle \langle 4 \rangle + \langle  1 \rangle \langle 3 \rangle \langle4 \rangle 
  + \langle 1\lt 3,2\lt3 \rangle \langle 4 \rangle - \langle 1,2 \rangle \langle 4  \rangle  \\ && -
  \langle 1\lt 4,2\lt 4 \rangle \langle  3\lt 4 \rangle + \langle 1,2 \rangle \langle 3  \rangle ; \\
   \partial (\langle 1 \rangle \langle 2, 3 \rangle \langle 4 \rangle ) &=& 
    \langle 2,3 \rangle \langle 4 \rangle - \langle 2,3 \rangle \langle 4 \rangle  -\langle 1\lt 2 \rangle \langle 3 \rangle \langle 4 \rangle 
    + \langle 1\rangle \langle 2\cdot 3 \rangle \langle 4 \rangle  -
    \langle 1 \rangle \langle 2 \rangle  \langle 4 \rangle \\ && - \langle 1 \lt 4 \rangle \langle 2 \lt 4, 3 \lt 4 \rangle  + \langle 1  \rangle \langle 2 , 3  \rangle ; \\ 
    \partial (\langle 1\rangle \langle 2 \rangle \langle 3,4 \rangle ) &=& 
    \langle 2\rangle \langle 3,4 \rangle - \langle 2\rangle \langle3,4 \rangle -\langle 1\lt 2 \rangle \langle 3,4 \rangle 
    + \langle 1\rangle \langle  3,4 \rangle  \\ &&+
   \langle 1\lt 3 \rangle \langle 2\lt 3 \rangle  
    \langle 4 \rangle -\langle 1 \rangle \langle 2 \rangle \langle 3 \cdot 4 \rangle + \langle 1 \rangle \langle 2 \rangle\langle 3 \rangle  ; \\
  \partial (\langle 1\rangle \langle 2 , 3,4 \rangle )
  & = &    \langle 2 , 3,4 \rangle -  \langle 2 , 3,4 \rangle - 
   \langle 1 \lt 2 \rangle  \langle 3,4 \rangle  +  \langle 1  \rangle \langle 2 \cdot 3,  4\rangle - \langle 1  \rangle \langle 2, 3 \cdot 4 \rangle + \langle 1  \rangle \langle 2,3 \rangle; \\
\partial (\langle 1\rangle \langle 2 \rangle \langle 3 \rangle \langle 4 \rangle ) &=& 
     \langle 2\rangle \langle 3 \rangle \langle 4 \rangle
      - 
     \langle 2\rangle \langle3 \rangle\langle 4 \rangle
      -\langle 1\lt 2 \rangle \langle 3 \rangle  \langle 4 \rangle
    + \langle 1\rangle \langle  3 \rangle \langle 4 \rangle
     +
    \langle 1\lt 3 \rangle \langle 2 \lt 3 \rangle \langle 4 \rangle -  \langle 1\rangle \langle 2 \rangle \langle 4 \rangle \\ &&
  -  \langle 1\lt 4 \rangle \langle 2 \lt 4  \rangle \langle 3 \lt 4  \rangle + \langle 1\rangle \langle 2 \rangle \langle 3 \rangle.
 \end{eqnarray*}
  
  \includegraphics[width=6.5in]{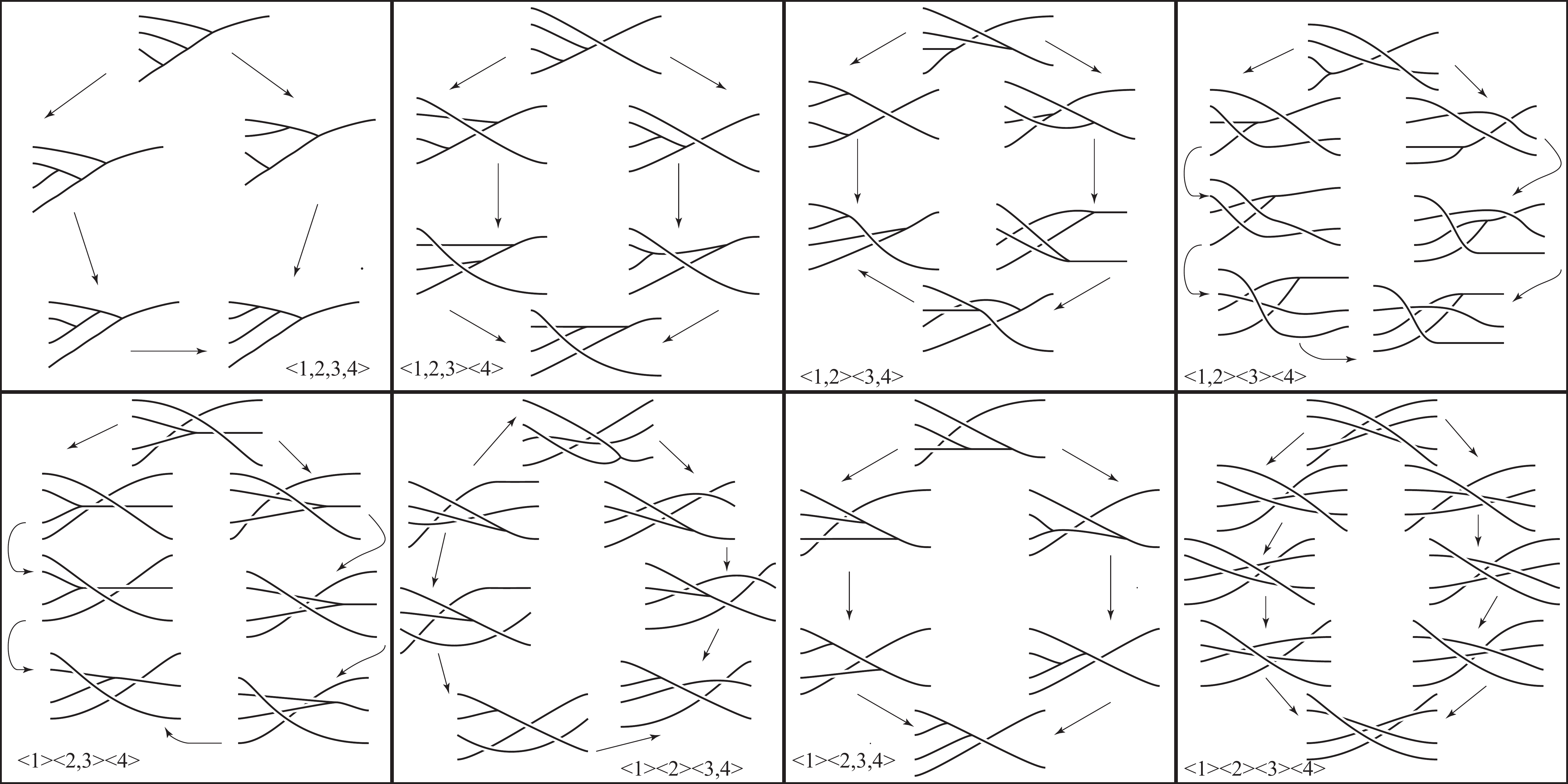}
  
  \vspace{1cm} 
  
  For any $a,b,c \in X$, we seek functions $\alpha: (\{a\}\times G)^3 \rightarrow A$,
  $\gamma_1: (\{a\}\times G)^2 \times  (\{b\}\times G)\rightarrow A$,
    $\gamma_2: (\{a\}\times G)\times  (\{b\}\times G)^2\rightarrow A$, and 
    $\theta: (\{a\} \times G) \times (\{b\} \times G) \times (\{c\} \times G) \rightarrow A$ that satisfy the following eight {\it cocycle conditions}:
   
  \begin{eqnarray*} \alpha(\langle 1 \cdot 2, 3,4 \rangle ) +  \alpha(\langle 1,2, 3 \cdot 4 \rangle )
  & = &
\alpha(  \langle 2,3,4 \rangle ) +\alpha( \langle 1, 2 \cdot 3,4 \rangle ) + \alpha(\langle 1,2,3 \rangle ) ; \\ 
 \gamma_1(\langle 1 \cdot 2, 3 \rangle \langle 4 \rangle ) + \gamma_1 ( \langle 1,2\rangle \langle 4 \rangle )
 + \alpha ( \langle 1\lt 4 , 2 \lt 4 , 3 \lt 4 \rangle )   &= & \gamma_1 ( \langle 2,3 \rangle \langle 4 \rangle ) +  \gamma_1( \langle 1, 2 \cdot 3 \rangle \langle 4 \rangle ) \\ && + \alpha( \langle 1 , 2  , 3 \rangle  ); \\ 
   \gamma_2( \langle 1 \cdot 2 \rangle \langle 3,4 \rangle )+
  \gamma_1(\langle 1,2 \rangle \langle 3 \cdot 4 \rangle )
   &=&
   \gamma_2(\langle 2 \rangle \langle 3,4 \rangle ) + \gamma_2( \langle  1 \rangle \langle 3,4 \rangle ) \\ &&
  + \gamma_1(\langle 1\lt 3,2\lt3 \rangle \langle 4 \rangle) +\gamma_1( \langle 1,2 \rangle \langle 3  \rangle ) ; \\  
  \theta( \langle 1 \cdot 2 \rangle \langle 3 \rangle \langle 4 \rangle) +  \gamma_1( \langle 1,2 \rangle \langle 4  \rangle )   +
\gamma_1(  \langle 1\lt 4,2\lt 4 \rangle \langle  3\lt 4 \rangle ) 
   &=& \theta( \langle 2 \rangle \langle 3 \rangle \langle 4 \rangle ) + \theta( \langle 1 \rangle \langle 3 \rangle \langle 4 \rangle ) \\ && 
   + \gamma_1(\langle 1\lt 3,2\lt3 \rangle \langle 4 \rangle) +\gamma_1( \langle 1,2 \rangle \langle 3  \rangle ) ;\\
   \theta(\langle 1\lt 2 \rangle \langle 3 \rangle \langle 4 \rangle ) +
   \theta( \langle 1 \rangle \langle 2 \rangle  \langle 4 \rangle )  +\gamma_2( \langle 1 \lt 4 \rangle \langle 2 \lt 4, 3 \lt 4 \rangle )
   &=& 
     \theta(\langle 1\rangle \langle 2\cdot 3 \rangle \langle 4 \rangle )  +\gamma_2( \langle 1  \rangle \langle 2 , 3  \rangle ); \\ 
  \theta( \langle 1\lt 3 \rangle \langle 2\lt 3 \rangle  
    \langle 4 \rangle )   +  \theta( \langle 1 \rangle \langle 2 \rangle\langle 3 \rangle  ) &=&
        \gamma_2( \langle 1\lt 2 \rangle \langle 3,4 \rangle )
       + \theta( \langle 1 \rangle \langle 2 \rangle \langle 3 \cdot 4 \rangle ) ; 
 \\ 
\gamma_2(  \langle 1 \lt 2 \rangle  \langle 3,4 \rangle  ) +\gamma_2( \langle 1  \rangle \langle 2, 3 \cdot 4 \rangle ) 
  & = &    
 \gamma_2( \langle 1  \rangle \langle 2 \cdot 3,  4\rangle )    + \gamma_2( \langle 1  \rangle \langle 2,3 \rangle ); \\
\theta( \langle 1\lt 2 \rangle \langle 3 \rangle  \langle 4 \rangle ) +  \theta(  \langle 1\rangle \langle 2 \rangle \langle 4 \rangle )   +  \theta( \langle 1\lt 4 \rangle \langle 2 \lt 4  \rangle \langle 3 \lt 4  \rangle  ) 
&=& 
    \theta( \langle 1\rangle \langle  3 \rangle \langle 4 \rangle )
     +
    \theta(  \langle 1\lt 3 \rangle \langle 2 \lt 3 \rangle \langle 4 \rangle )  \\ &&
+  \theta( \langle 1\rangle \langle 2 \rangle \langle 3 \rangle ) . \end{eqnarray*} 

We note that these cocycle conditions are given in general, and the normal orientations in the figures are chosen to be pointing upwards. In the case that the coloring is by $\tilde{R}$, the normals to sheets coincident to arcs have to be chosen consistently. In this case, the signs on some of the triple points reverse. Thus the cocycle conditions that are below reflect these changes in signs. 
  
 In the general theory of knotted foams, there are isotopy moves that are analogous to (and include) the Roseman moves. When an embedded foam is projected to $3$-space its $0$-dimensional multiple points will consist of the four scenarios that are depicted above, branch points induced by Reidemeister type-I moves and the isolated vertex that occurs when a trivalent vertex undergoes a twist. These last two singularities are depicted here.  

\includegraphics[width=3in]{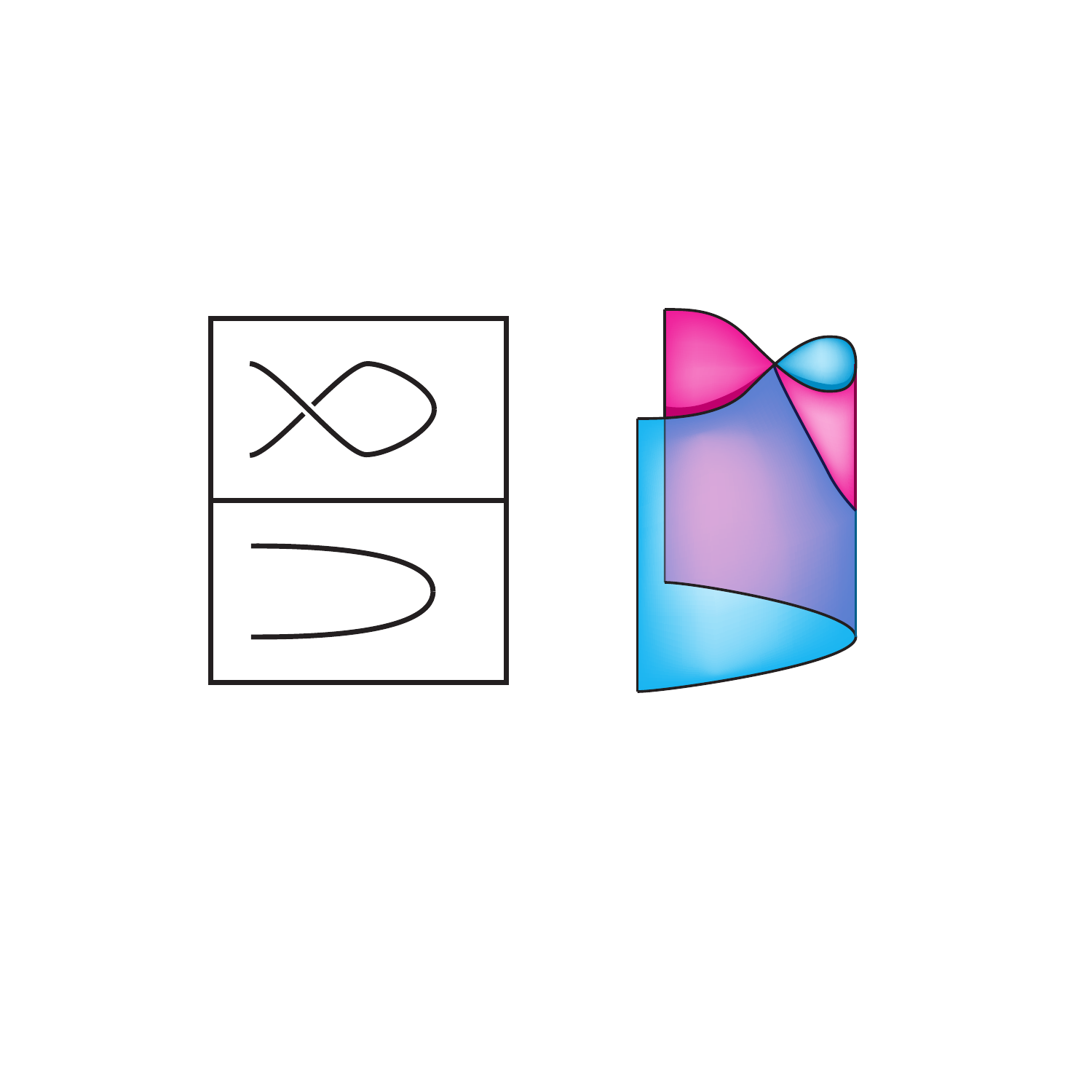}
\includegraphics[width=3in]{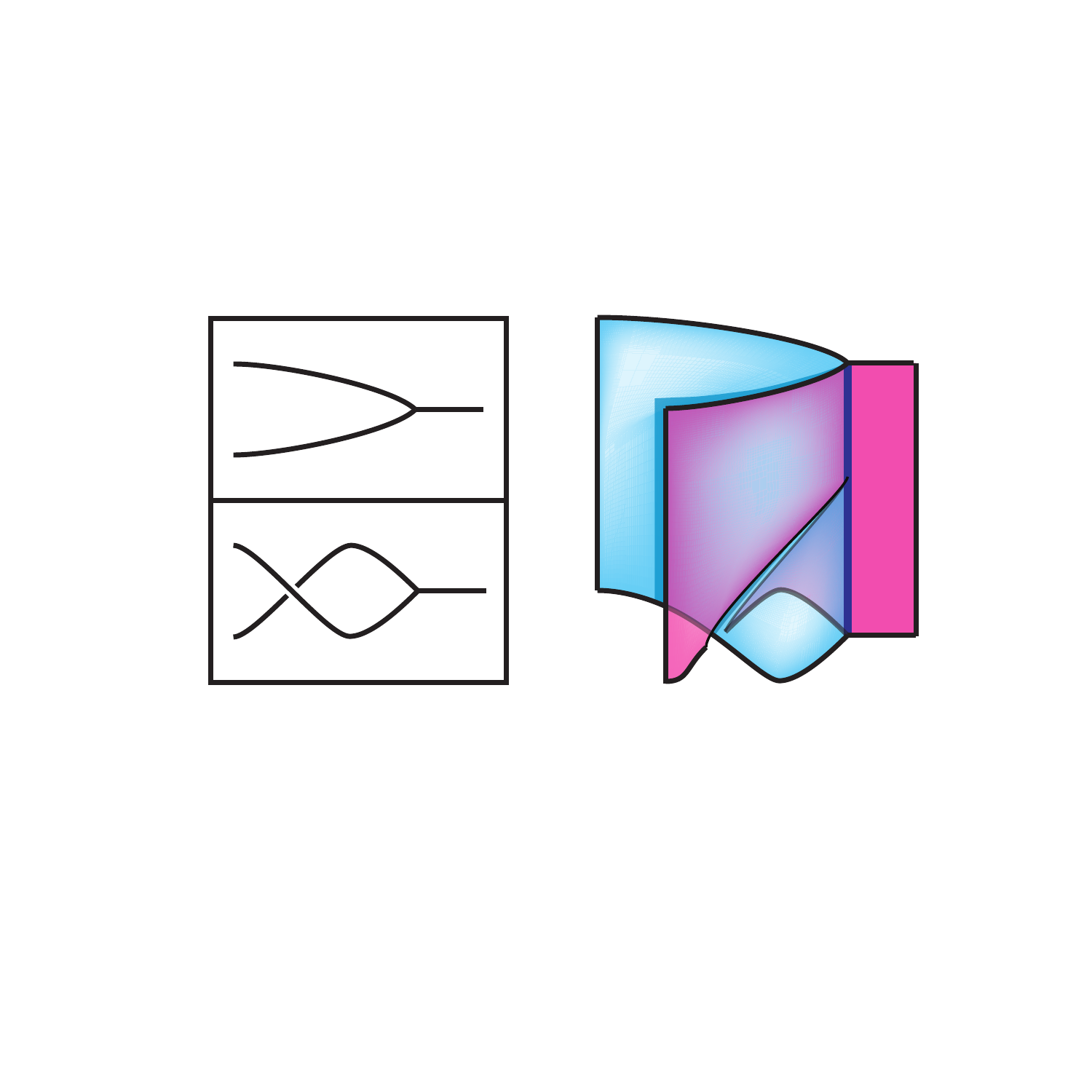}

\vspace{-2cm}

Just as a trivalent graph represents an embedded handlebody in $3$-space, the embedded foams represent certain $4$-manifolds with boundary that are embedded in $4$-space. The $4$-manifolds are regular neighborhoods of the foams. As such, the Matveev-Piergallini \cite{Mat2,Piergallini} moves for special spines can be applied to the underlying foams without changing the $4$-manifold. These moves are the  move $\langle 1,2,3,4 \rangle$ indicated above and the invertibility condition for the basic foam $Y^2$.
This condition and other obvious invertibilities are indicated here. 

\begin{center}
\includegraphics[width=4.5in]{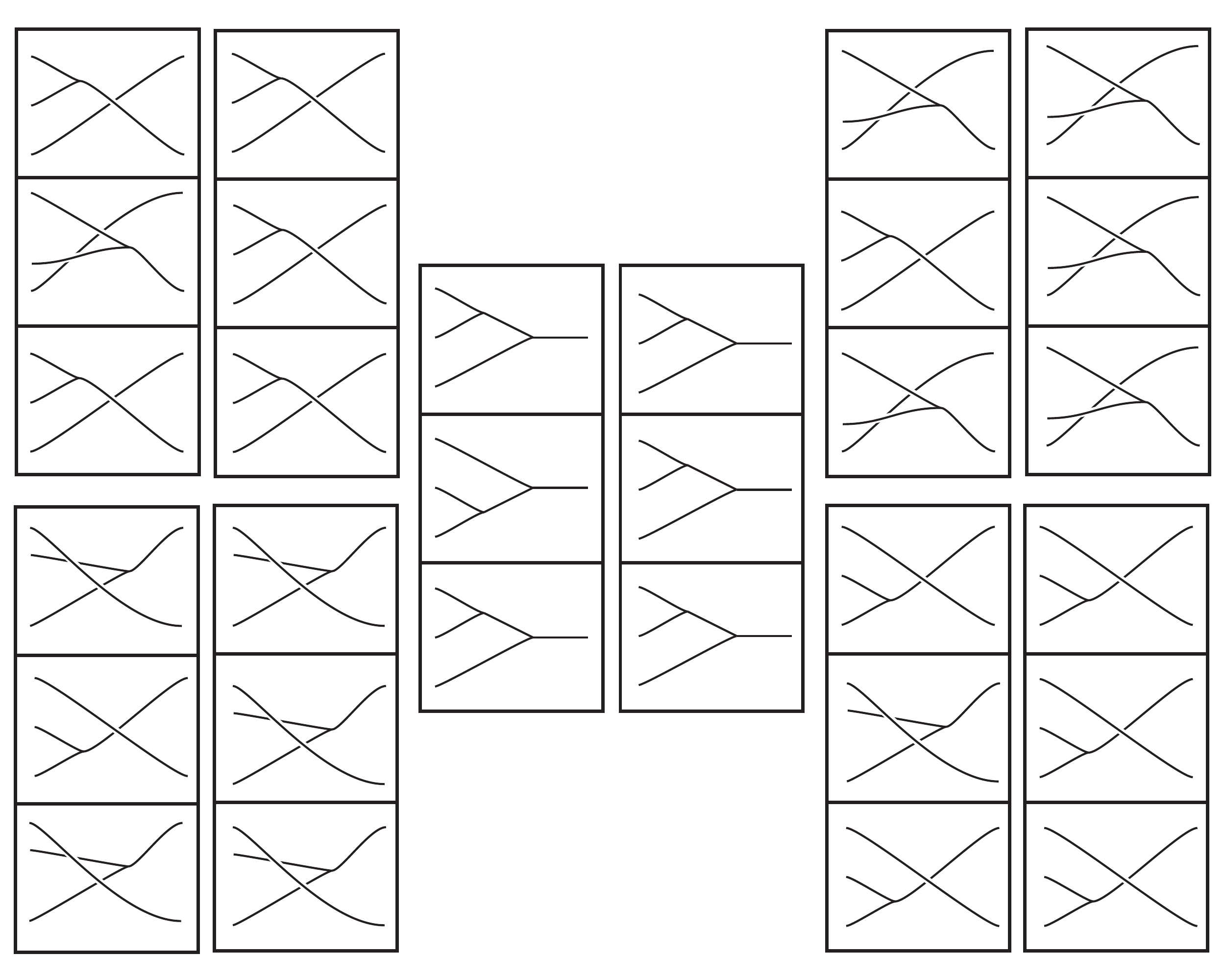}
\end{center}

The analogues of the Roseman moves, then, are (1) the invertibility of each $0$-dimensional multiple point including those depicted and the invertibility of the twisted vertex (both elliptically and hyperbolically), (2) the eight movie-moves indicated above, (3) The original seven Roseman moves, and (4) pushing a twisted vertex through a transverse sheet. The proof that these moves suffice will be presented elsewhere. 

\section{Cocycle Invariants of knotted foams}
The cocycle conditions for the associated quandle for a $G$-family of quandles give us quantities that are invariant under the eight main movie moves. Here we are considering labeling the arcs on the left of the string diagrams from bottom to top with elements of the associated quandle of a $G$-family of quandles. When arcs are conjoined by trivalent vertices, the 
$X$-coloring is monochromatic while the 
$G$-coloring varies. Thus quandle cocycle invariants can be defined for knotted foams in the following way: 
\begin{itemize}
\item choose a coloring of the foam by a $G$-family of quandles;
\item assign cocycle values at the $0$-dimensional multiple points of the foam --- these are points at which a $Y^1$ crosses a transverse sheet, a vertex of the dual to the tetrahedron, or a triple point;
\item take the product of the cocycle values over all the $0$-dimensional multiple points of the closed foam;
\item take the multi-set of values over all colorings.
\end{itemize}
That this process is invariant depends upon the cocycle conditions, the assignation of signs to the $0$-dimensional multiple points, the existence of a good involution \cite{IIJO}
on a $G$-family, and the vanishing of the cocycles upon degenerate chains. 

For the purposes of the current paper, we will assume that $\alpha$, $\gamma_1$, and $\gamma_2$ all are constantly and trivially valued to be $0$. 
With the convention that the normal-orientation follows consistently the two of three sheets that are labeled by $(a,1)\in \tilde{R}$ at the junction of three sheets, the cocycle conditions read as follows:

$\langle 1,2 \rangle \langle 3\rangle \langle 4\rangle:$
 \begin{eqnarray*}
  \theta( (a,1),(c,k), (d,\ell)) 
   &=& \theta( (a,0),(c,k),(d,\ell) ) + \theta( (a,1),(c,k),(d,\ell)); \\   
     \theta( (a,0),(c,k), (d,\ell)) 
   &=& \theta( (a,1),(c,k),(d,\ell) ) - \theta( (a,1),(c,k),(d,\ell)).  \end{eqnarray*} 
   
   $\langle 1 \rangle \langle2, 3\rangle \langle 4\rangle:$
   \begin{eqnarray*}
  - \theta((a, g), (b,1),(d,\ell))
    +
   \theta( (a,g),(b,1),(d,\ell))   
    &=& 
     \theta((a,g),(b,0),(d,\ell)); \\ 
       \theta((2b-a, g), (b,1),(d,\ell))
    -
   \theta( (2b-a,g),(b,1),(d,\ell))   
    &=& 
     \theta((a,g),(b,0),(d,\ell)); \\
    \theta((a, g), (b,1),(d,\ell))
    +
   \theta( (a,g),(b,0),(d,\ell))   
    &=& 
     \theta((a,g),(b,1),(d,\ell)); \\  
      \theta((2b-a, g), (b,0),(d,\ell))
    +
   \theta( (a,g),(b,1),(d,\ell)).   
    &=& 
     \theta((a,g),(b,1),(d,\ell) \end{eqnarray*} 
     
 $\langle 1 \rangle \langle2\rangle \langle 3, 4\rangle:$    
 \begin{eqnarray*}    
    \theta((a,g),(b,h),(c,1))  &=&
      \theta((a,g),(b,h),(c,0) ) +    \theta((2c-a,g),(2c-b,h),(c,0) ); 
 \\ 
 \theta((a,g),(b,h),(c,1))  &=&
      \theta((a,g),(b,h),(c,0) ) +    \theta((a,g),(b,h),(c,1) ); 
 \\ 
 \theta((a,g),(b,h),(c,0))  &=&
      \theta((a,g),(b,h),(c,1) ) -   \theta((a,g),(b,h),(c,1) ); 
 \\ 
  \theta((a,g),(b,h),(c,0))  &=&
     - \theta((2c-a,g),(2c-b,h),(c,1) ) +   \theta((2c-a,g),(2c-b,h),(c,1) ). \end{eqnarray*}
 
  $\langle 1 \rangle \langle2\rangle \langle 3\rangle \langle4\rangle:$ 
 \begin{eqnarray*}
\theta( \langle 1\lt 2 \rangle \langle 3 \rangle  \langle 4 \rangle ) +  \theta(  \langle 1\rangle \langle 2 \rangle \langle 4 \rangle )   +  \theta( \langle 1\lt 4 \rangle \langle 2 \lt 4  \rangle \langle 3 \lt 4  \rangle  ) 
&=& 
    \theta( \langle 1\rangle \langle  3 \rangle \langle 4 \rangle )
     +
    \theta(  \langle 1\lt 3 \rangle \langle 2 \lt 3 \rangle \langle 4 \rangle )  \\ &&
+  \theta( \langle 1\rangle \langle 2 \rangle \langle 3 \rangle ) . \end{eqnarray*}

Consider the case, in which the group 
$G=\mathbb{Z}_2=\{0,1\}$, the set $X=R_3=\{0,1,2\}$, and the associated quandle is $Q=X\times G$. The quandle actions on $R_3$ are 
$a \lt_0 b =a$, $a\lt_1 b=2b-a$, and the action on $Q$ is $(a,g)\lt(b,h)=(a\lt_h b, g)$. Consider further, 
Mochizuki's $3$-cocycle $\theta_p:X^3\to\mathbb{Z}_3$ which, in this case, can be simplified to 
$$ \theta_3(a,b,c):=(a-b)(c^3 +c^2b +  b^2c). \]
Then
\[ \theta((x_1,g_1),(x_2,g_2),(x_3,g_3))
= \left\{ \begin{array}{lr}
\theta_3(x_1,x_2,x_3) & \ \ {\mbox{\rm if}} \ \ g_1=g_2=g_3=1, \\
0 & \ \ {\mbox{\rm otherwise.}}\end{array}\right.
 $$

\begin{lemma} The function $\theta$ satisfies the cocycle conditions above. \end{lemma}
{\sc Proof.} The cocycle conditions  for the relations 
$\langle 1,2 \rangle \langle 3\rangle \langle 4\rangle,$  $\langle 1 \rangle \langle 2,3\rangle \langle 4\rangle,$ and $\langle 1 \rangle \langle 2\rangle \langle 3,4\rangle,$ all follow easily since all expressions that involve $0$ as a group element are trivial. The last condition, $\langle 1 \rangle \langle 2\rangle \langle 3\rangle \langle 4\rangle,$  involves trivial terms if any one of the arguments has a $0$ as the group element. Otherwise, all the group elements are $1$, and the result follows because Mochizuki's function satisfies this particular cocycle condition. $\Box$

For the reader's convenience, we tabulate values of the cocycle $\theta$. Since $\theta=0$ if any of $g,h,k=0$ only values for which $g=h=k=1$ are indicated. Similarly, when $x=y$ or $y=z$, these values are excluded.

\begin{center}
\begin{tabular}{||c|c||} \hline \hline 
$((x,g),(y,h),(z,k))$ & $\theta((x,g),(y,h),(z,k))$ \\ \hline \hline 
$((0,1),(1,1),(2,1))$ & $1$       \\ \hline 
$((0,1),(2,1),(1,1))$ & $1$      \\ \hline 
$((1,1),(0,1),(2,1))$ & $2$       \\ \hline 
$((1,1),(2,1),(0,1))$ & $0$      \\ \hline 
$((2,1),(0,1),(1,1))$ & $2$       \\ \hline 
$((2,1),(1,1),(0,1))$ & $0$       \\ \hline 
$((0,1),(1,1),(0,1))$ & $0$      \\ \hline 
$((0,1),(2,1),(0,1))$ & $0$     \\ \hline 
$((1,1),(0,1),(1,1))$ & $1$     \\ \hline 
$((1,1),(2,1),(1,1))$ & $2$    \\ \hline 
$((2,1),(0,1),(2,1))$ & $1$   \\ \hline 
$((2,1),(1,1),(2,1))$ & $2$    \\ \hline \hline
\end{tabular}
\end{center}

\section{The value of the invariant}

Our initial example of a knotted foam is the 
$2$-twist spin of the knotted trivalent graph that represents the knotted handlebody $5_2$ in the tables~\cite{IKMS}. Here we demonstrate that the cocycle value is non-trivial by exhibiting the double decker set in the presence of a non-trivial coloring. The coloring $(R,B,G)$ represents any coloring for which all three colors are different. It is not difficult to observe that if any two of these colors are coincident, then they all are the same. 

The $2$-twist-spin of the foam is illustrated in a step-by-step process. From this ``movie", we can construct the relevant part of the decker-set --- the preimage of the double points on the abstract foam. The coloring chosen is given so that the top arc is colored $(z,0)$. (In the figure, we name this arc the ``north arc"). Consequently, any triple points that involve this arc will not contribute to the cocycle invariant. The ``south arc" winds from the bottom left to the bottom right of the figure. The ``east arc'' is on the right of the figure. The central arc is the remaining arc. On the decker set, the lower decker set that involves the north arc is indicated with a thin line. At the pre-image of the triple points on the lowest sheet, the source region is indicated with a black dot. Sign conventions follow those in \cite{CS:Book}.  

\begin{center}
\includegraphics[width=6.5in]{nuspinmovie}
\end{center}

\begin{center}
\includegraphics[width=6.0in]{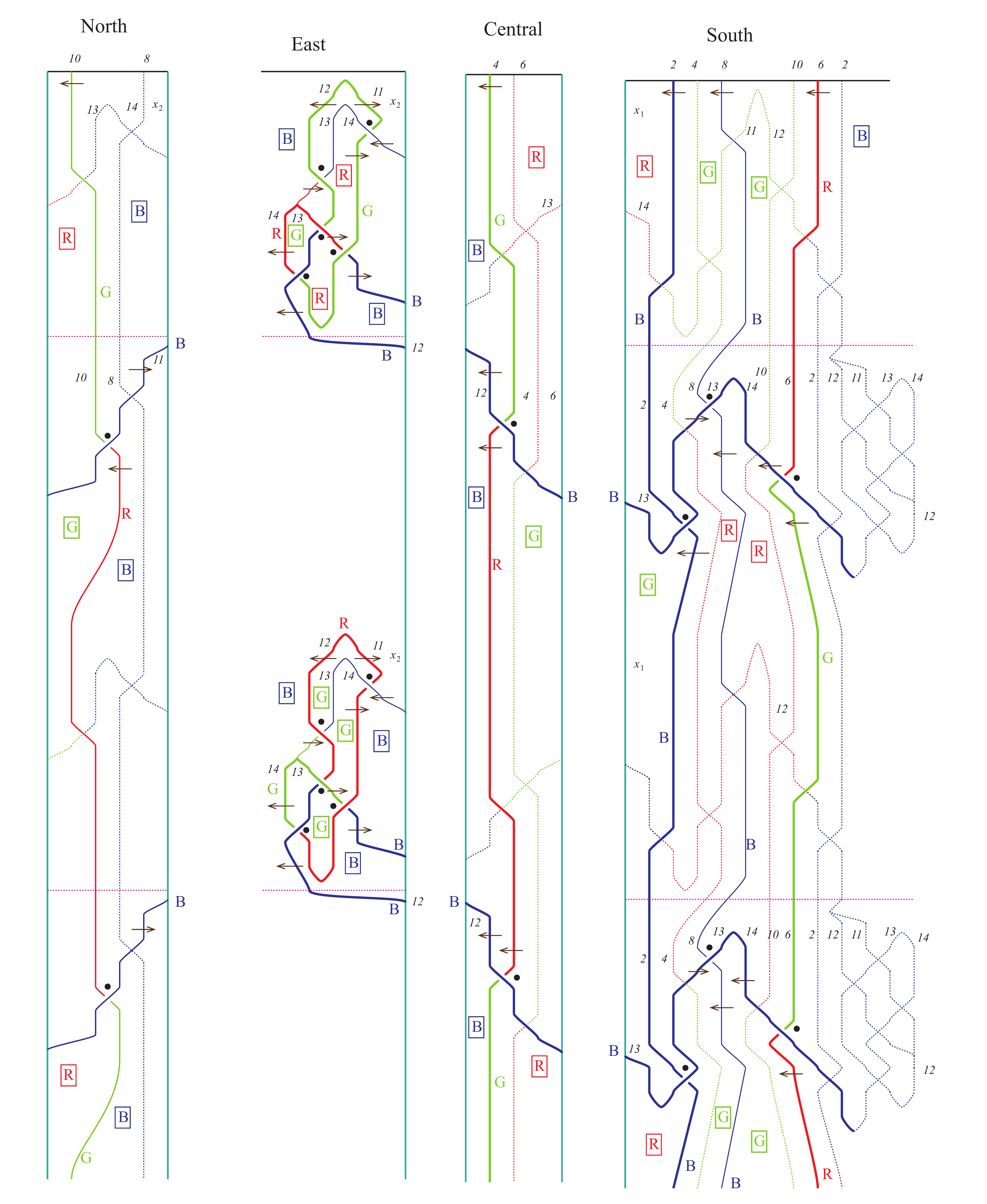}
\end{center}

From the decker-set, we obtain the value
$$ +\theta((x,1),(z,1),(y,0))
+\theta((x,1),(y,0),(z,1))
+\theta((x,1),(y,1),(x,1))
+\theta((x,1),(x,1),(z,1))$$
$$+\theta((x,1),(z,1),(y,1))-\theta((z,1),(y,0),(y,1))
-\theta((x,1),(z,1),(y,1))
-\theta((y,0),(z,0),(y,0))$$
$$-\theta((y,1),(x,1),(y,1))
-\theta((z,1),(y,1),(y,1))
 +\theta((z,1),(x,1),(y,0))
+\theta((z,1),(y,0),(x,1))$$
$$+\theta((z,1),(y,1),(z,1))
+\theta((z,1),(z,1),(x,1))
+\theta((z,1),(x,1),(y,1))-\theta((x,1),(y,0),(y,1))$$
$$-\theta((z,1),(x,1),(y,1))
-\theta((y,0),(x,0),(y,0))
-\theta((y,1),(z,1),(y,1))
-\theta((x,1),(y,1),(y,1))$$

$$= 
\theta((x,1),(y,1),(x,1))
-\theta((y,1),(x,1),(y,1))+\theta((z,1),(y,1),(z,1))
-\theta((y,1),(z,1),(y,1)).
$$
 
 By evaluating this sum for all possible $x,y,z$ with these values distinct, we obtain the value $2$ for each non-trivial color. We conclude  with the following result.

 \begin{center}
\includegraphics[width=4in]{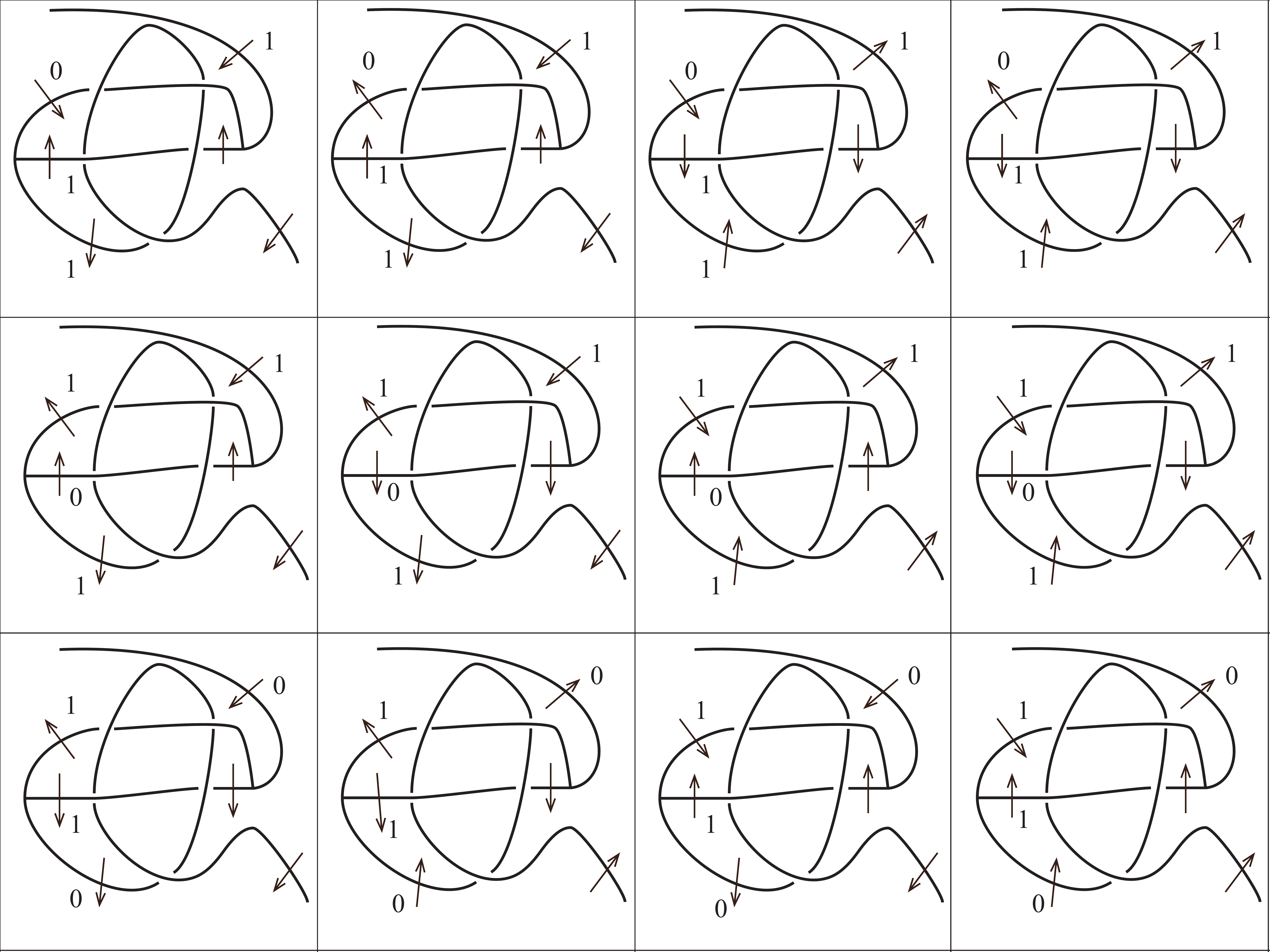}
\end{center}
 
 \begin{theorem} The quandle cocycle value for the knotted foam illustrated is given by 
$60+12t+12t^2$ where we indicate the multiset of values by means of a polynomial. \end{theorem}
{\sc Proof.} There are 12 possible flows that are compatible for colorings by $\tilde{R} = \{ (0,0),$ $(0,1),$ $(1,0),$ $(1,1),$ $(2,0),$ $(2,1)\}$ with quandle rules
$(a,g)\lt(b,0)=(a,g)$ and $(a,g)\lt(b,1)=(2b-a,g)$
for $a,b\in\Z_3$ in which two edges have group elements 1 upon them. These are illustrated in the figure above the statement of the theorem. For any flow in the bottom two rows, the colorings of the arcs must be trivial. There are 24 such trivial colors. 

In the top row, there are $3\times 4=12$ trivial colorings. 

In addition, there are 24 trivial colors from all three edges being  colored with $(a,0)$  for $a\in \{0,1,2\}$ and eight different possible transverse orientations.

We have indicated --- by means of the decker set and $(R,B,G)$ --- the cocycle invariants of three of the colorings associated to the orientation in top left corner of the table. Each gives a value of $2\in \Z_3$. The $(1,2)$ entry will give the same values. The last two rows will result in all of the triple points having their orientations reversed, and give the other 12 non-trivial values of $1$. $\Box$

\section{Summary}

We have illustrated a knotted foam that has non-trivial cocycle invariant. In the immediate future, a proof that the foam moves summarized here are a complete set of Reidemeister/Roseman moves for knotted foams in $4$-space will be presented. Further examples of knotted foams with non-trivial cocycle invariants will be given. Finally, a serious investigation on the relationships between group and quandle cohomology is due.

\begin{flushleft}
J. Scott Carter \\ 
Department of Mathematics \\ 
University of South Alabama \\ 
Mobile, AL 36688 
USA \\
E-mail address: {\tt carter@southalabama.edu} 
\end{flushleft}

\begin{flushleft}
Atsushi Ishii  \\University of Tsukuba\\ 1-1-1 Tennodai \\ Tsukuba, Ibaraki 305-8571, Japan\\ 
E-mail address: {\tt  aishii@math.tsukuba.ac.jp }
\end{flushleft}
\end{document}